\documentclass{article}
\usepackage{latexsym}
\usepackage{amsfonts,amsmath,amsthm}
\usepackage{amssymb}
\usepackage{amscd}

\makeatletter
\def\theequation{\arabic{section}.\arabic{equation}}
\@addtoreset{equation}{section}
\newtheorem{theorem}{Theorem}[section]

\newtheorem{proposition}[theorem]{Proposition}

\theoremstyle{definition}
\newtheorem{definition}[theorem]{Definition}
\theoremstyle{remark}

\newtheorem{remark}[theorem]{Remark}

\newcommand{\de}{\mathrm{d}} 

\newcommand{\be}{\begin{equation}}
\newcommand{\ee}{\end{equation}}
\newcommand{\R}{\mathbb{R}}
\newcommand{\N}{\mathbb{N}}
\newcommand{\inte}{\int_{0}^{1}}

\newcommand{\tig}{\tilde{g}_0}
\newcommand{\M}{{\mathcal M}}
\newcommand{\J}{{\mathcal J}}

\begin{document}

%-----------------------------------------------------------------

\title{Normal geodesics connecting\\ two non--necessarily 
spacelike submanifolds\\ in a stationary spacetime
\footnote{Work supported by M.I.U.R. Research Project PRIN07 ``Metodi Variazionali e 
Topologici nello Studio di Fenomeni Nonlineari''.}}
\author{Rossella Bartolo$^{\ddag,1}$, Anna Maria Candela$^{\dag,2}$,
Erasmo Caponio$^{\ddag,3}$\\
{\it\small $^\ddag$Dipartimento di Matematica, Politecnico di Bari}, \\
{\it\small Via E. Orabona 4, 70125 Bari, Italy}\\
{\it\small $^\dag$Dipartimento di Matematica, Universit\`a degli Studi di Bari}, \\
{\it\small Via E. Orabona 4, 70125 Bari, Italy}\\
{\it\small e-mails: $^1$r.bartolo@poliba.it, $^2$candela@dm.uniba.it, $^3$caponio@poliba.it}
\vspace{1mm}}

\date{}
\maketitle

\begin{center}
{\bf\small Abstract}
\vspace{3mm}
\hspace{.05in}\parbox{4.5in}{{\small In this paper we obtain an existence theorem for normal geodesics
joining two given submanifolds in a globally
hyperbolic stationary spacetime $\M$. The proof is based on both
variational and geometric arguments involving the causal structure of $\M$,
the completeness of suitable Finsler metrics associated to it and some basic properties of a submersion. By this interaction, unlike
previous results on the topic, also non--spacelike submanifolds can be handled.}}
\end{center}

\noindent
{\it \footnotesize Key words.} {\scriptsize Stationary spacetime, 
normal geodesic, global hyperbolicity,  
Finsler metric, submersion. }\\
{\it \footnotesize 2000 MSC.} {\scriptsize 53C50, 53C22, 58E10, 53C60}.

%-----------------------------------------------------------------

\section{\bf Introduction and background tools}\label{s1}

The aim of this paper is  proving an existence result for normal geodesics
connecting two rather general submanifolds in a globally hyperbolic stationary spacetime.
On one side, in the case when both the submanifolds are compact, our result is independent of
anyone of the global splittings that  such a type of spacetime admits
for the global hyperbolicity assumption.
On the other hand we are able to treat, for particular splittings,
also the case when one of the submanifolds is non--compact and non--spacelike.

Let us recall the basic notions related to this topic
(cf. \cite{bee} for the background material on Lorentzian geometry used throughout the paper).

A {\em Lorentzian manifold} $({\mathcal M},g)$ is a smooth
connected  finite dimensional mani\-fold  equipped with a $(0,2)$ symmetric
non--degenerate tensor field $g$ having index $1$.

A {\em geodesic} of $({\mathcal M},g)$ is a smooth curve
$z: [a,b]\subset\R\to \M$ satisfying the equation
\[
\nabla_{s}\dot z\ =\ 0,
\]
where $\nabla_s$ is the covariant derivative along $z$ associated to the Levi--Civita
connection of the metric $g$.
Without loss of generality, we can reduce our study
to geodesics defined in the same interval $[0,1]$; furthermore,
it is well known that geodesics satisfy the conservation law
$g(z)[\dot z,\dot z] = E_{z}$.
Thus, they are classified according to their {\em causal character},
that is according to the sign of the constant $E_{z}$: $z$ is said {\em
timelike} if $E_{z}<0$, {\em lightlike} if $E_{z}=0$, {\em
spacelike} if $E_{z}>0$ or $\dot z= 0$, {\em causal} if $E_z\leq 0$.
The same terminology is used also for any vector and
for any vector field if it has the same causal character at each point,
for any piecewise smooth curve (according to the causal character
of its velocity vector field) and for submanifolds.
In particular, a submanifold $P$ of $\mathcal M$
is {\em spacelike} if $g$ restricted to $T_pP$ is positive definite for each $p\in P$.

A {\em spacetime} is a Lorentzian manifold with a prescribed time--orien\-ta\-tion,
that is with a continuous choice of a causal cone at each point of $\mathcal M$.
In such a case a piecewise smooth causal curve on $\mathcal M$ is said
{\em future--pointing} (resp. {\em past--pointing})
 if its velocity vector field belongs
to the cones labeled as {\em future} ones
at any point where it is defined.

A vector field $K$ on $\mathcal M$ is {\em Killing} if
one of the following equivalent assertions holds true (see \cite[Propositions 9.23 and 9.25]{o}):
\begin{itemize}
\item[{\sl (i)}] the stages of its local flow consist of isometries;
\item[{\sl (ii)}] the Lie derivative of $g$ in its direction is $0$;
\item[{\sl (iii)}] $g[\nabla_X K,Y] = - g[\nabla_Y K,X]$ for each pair of vector fields $X,Y$.
\end{itemize}
It is easy to see that if $K$ is a Killing vector field and
$z$ is a geodesic, then a constant $C_z\in\R$ exists such that
\be\label{vinc}
g(z)[\dot z,K]= C_z.
\ee
The existence of a timelike Killing vector field
gives some important information on the structure of the
manifold (e.g., cf. \cite{s1}). Moreover, observers traveling
 on integral curves of timelike Killing vector fields
see a constant metric.

A spacetime $(\mathcal M,g)$ is called {\em stationary} if it admits a timelike Killing
vector field. It is {\em globally hyperbolic} if it admits a (smooth) spacelike Cauchy hypersurface,
i.e. a subset crossed exactly once by any inextensible timelike curve.

If a globally hyperbolic stationary spacetime $\mathcal M$ admits
at least one complete Killing vector field $K$ (i.e., the integral curves of $K$
are defined on $\R$), then it is {\em standard stationary} (see \cite[Theorem 2.3]{cfs}),
that is $\mathcal M$ splits as a product $\M_0\times\R$, where the connected finite dimensional
manifold $\M_0$ is endowed with a Riemannian metric $g_0$ and the metric $g$ is given by
\be\label{stand}
g(x,t)[(y,\tau),(y,\tau)]=g_0(x)[y,y]+ 2g_0(x)[\delta(x),y]\tau-\beta(x)\tau^2
\ee
for any $(x,t)\in \M_0\times \R$, $(y,\tau)\in T_x\M_0\times\R$,
with $\delta$ vector field and $\beta$ positive function, both on $\M_0$;
in this case $K=\partial_t$. It is well known that
locally any stationary spacetime looks like a standard one.\footnote{See \cite{JavS08} for a
characterization in terms of the causal properties that a stationary Lorentzian manifold with
a complete timelike Killing vector field has to satisfy in order to be standard stationary.}

In \cite{cfs} it is proved that a globally hyperbolic stationary spacetime, endowed with a complete timelike
Killing vector field $K$ (and, thus, standard stationary) and with a complete
spacelike Cauchy hypersurface $S$, is geodesically connected. Remarkably enough, there are
counterexamples to geodesic connectedness if one of the assumptions in \cite{cfs} is dropped.

Variational methods are already used in \cite{gm} for studying the geodesic
connectedness in standard stationary spacetimes, possibly with boundary.
In that paper, the authors introduce a variational principle for geodesics based
on the natural constraint \eqref{vinc} and prove the geodesic connectedness
with respect to the metric \eqref{stand}
under boundedness  assumptions on the vector field $\delta$ and on the scalar field $\beta$:
\be\label{GM}
g_0(x)[\delta(x),\delta(x)]\leq C,\quad\quad\quad m_1\leq \beta(x)\leq m_2,\ee
 for some $C$, $m_1$, $m_2 > 0$ and for any $x\in \mathcal M_0$.
Obviously, the hypothesis \eqref{GM}, and then the result in \cite{gm}, depends on the given global splitting $\mathcal M_0\times\R$ of
the stationary spacetime $\mathcal M$.

In the subsequent paper \cite{gp} an {\em intrinsic} approach to the problem of geodes\-ics connectedness is developed.
Namely, the variational principle in \cite{gm} is translated in a splitting independent form
(similar to Theorem~\ref{t1} but for the case of the two--point boundary conditions)
and a compactness assumption on the infinite dimensional manifold of the paths between two points,
called {\em pseudocoercivity}, is introduced. Such assumption implies global hyperbolicity, but it is
rather difficult to establish if it holds and, indeed, in order to furnish an example of a stationary
spacetime satisfying pseudocoercivity, a standard stationary one is  chosen.

In \cite{cfs} the authors essentially show that their intrinsic geometric assumptions, involving the causal structure
of the spacetime, are equivalent to pseudocoercivity. For a given complete spacelike smooth Cauchy
hypersurface $S$ they consider the  manifold $S\times\R$ which, by the flow of $K$, is diffeomorphic
to $\M$ and isometric to $(\M,g)$ if endowed with a metric as in (\ref{stand}) such that
\be\label{iden}
\begin{split}
&S\equiv\M_0,\quad K\equiv \partial_t,\\
&x=\pi_S(z)\quad \hbox{with $\pi_S$ canonical projection on $S$},\\
&g(z)[K(z),K(z)] = - \beta(x),\\
&\delta(x)\; \text{is the orthogonal projection of $K(z)$ on $T_zS$, for any $z\in S$.}
\end{split}
\ee
Even if, in general, the global splitting is not unique and not
 canonically associated to $\M$,
the result obtained is independent of the chosen $K$ and $S$ and
no growth hypothesis on the coefficients of the metric needs to be assumed.

Here, our aim is stating an existence result for
{\em normal geodesics} joining two fixed submanifolds $P$ and $Q$,
i.e. geodesics $z:[0,1]\rightarrow \mathcal M$ such that
\be\label{orto}
\begin{cases}
z(0)\in P,\ z(1)\in Q,\\
\dot z(0)\in T_{z(0)}P^\perp,\ \dot z(1)\in T_{z(1)}Q^\perp .
\end{cases}
\ee
We assume that the spacetime $\M$ satisfies
the assumptions in \cite{cfs}, that is $M$ is a stationary Lorentzian manifold
endowed with a complete timelike Killing vector field $K$ and a smooth spacelike
Cauchy hypersurface $S$. Then, following \cite{cym},
to each of the global splittings that $\mathcal M$ admits we can associate
two Finsler metrics of Randers type,
named {\em Fermat metrics} (see Section \ref{secfermat}).
Such metrics are related to the Fermat principle for future--pointing and past--pointing
lightlike geodesics  and
their  completeness is linked to the global hyperbolicity of $\M$ (see \cite[Theorem 4.8]{cym}).

We also remark that a standard stationary spacetime can  be seen as the total space
of a Lorentzian submersion $\pi\colon \M \to  \M_0$ where the  one--dimensional fibers
are the flow lines of $\partial_t$ (see \cite[Example 3.2]{CaJaPi09}).

The use of the Fermat metrics and the properties of a submersion seem to be very convenient for handling
our problem  as the causal techniques used in \cite{cfs} seem not easily extensible from
the two--point boundary conditions to the boundary data (\ref{orto}) (cf. the proof of \cite[Lemma 5.5]{cfs}).

It is worth to stress that in Theorem \ref{tm} below we deal with rather general submanifolds
which do not appear in related papers on the topic.

Denoting by $\Psi\colon\R\times \M\to\M$ the flow of $K$  and considering
a smooth submanifold $R$ of $\mathcal M$,
we can introduce the continuous function $s_R\colon R\to \R$ defined as follows:
for each $r \in R$, let $s_R(r)$ be the value of the parameter of the flow of $K$ such that
\[
\{\Psi(s_R(r),r)\}= \Psi(\mathbb R\times\{r\})\cap S.
\]
Observe that $s_R$ is well defined and continuous since the flow lines of $K$, being timelike curves, intersect
the Cauchy hypersurface $S$ in a unique point.

We assume that $P$ and $Q$ are two smooth immersed submanifolds
which are disjoint connected and closed as topological subspaces of $\mathcal M$ and
which satisfy one of the following conditions:
\begin{itemize}
\item[$(H_1)$] $P$ is compact (as a topological subspace of $\mathcal M$)
and
\be\label{sQ}
\sup_{q\in Q}|s_Q(q)|=D_Q<+\infty;
\ee
\item[$(H_2)$] two smooth submanifolds $P_S$ and $Q_S$ of $S$ exist, such
that one of them is compact (as a topological subspace of  $\mathcal{M}$)
and
\[
P=\Psi(\R\times P_S), \qquad Q=\Psi(\R\times Q_S).
\]
\end{itemize}

Notice that if in $(H_1)$ $Q$  is a compact submanifold as well,
assumption \eqref{sQ} is satisfied for any Cauchy hypersurface $S$.
Moreover, as $P$ is compact we have that
\be\label{proiezionecompatta}
\Psi(G_P)\subset S \quad \hbox{is a compact subset of $\mathcal M$,}
\ee
where $G_P = \{(s_P(p),p): \ p \in P\}$ is the graph of $s_P$.

In the previous literature on this subject, which is also mainly
concerned with a fixed a priori splitting, much more restrictive
assumptions  are  imposed on the submanifolds $P$ and $Q$ in order
to apply variational methods; such assumptions  make impossible to
handle submanifolds as in hypotheses $(H_1)$ or $(H_2)$.

 For example, in \cite{c} it is considered a standard {\em static}
spacetime (i.e., a standard stationary one with $\delta=0$) and the
submanifolds $P$ and $Q$ are given as $P=S_1\times \{t_p\}$,
$Q=S_2\times \{t_q\}$, with $S_1$, $S_2$ submanifolds of $ \M_0$,
$t_p,t_q\in\R$. Moreover, some results  are  obtained in standard
stationary spacetimes if again  $P=S_1\times \{t_p\}$ and
$Q=S_2\times\R$ (see \cite{cs} for  lightlike geodesics, \cite{cs1}
for spacelike  ones  and \cite{cms} for a result in the orthogonal
splitting case). As in the seminal paper by K. Grove (in the
Riemannian setting) \cite{g}, the submanifolds $S_1,S_2$ are always
assumed to be closed and at least one of them has to be compact,
although it is also possible to consider more general cases, up to
suitable additional assumptions involving them (cf.  \cite{bgs1} and
references therein).

Following the ideas developed in \cite{gp}, in \cite{bgs} is stated
a result for geodesics joining {\em spacelike} submanifolds of a
stationary spacetime, and again in the standard  case (cf.
\cite[Appendix B]{bgs}) such submanifolds turn out to be as in the
above cited papers.

Now, we are ready to state our main result.

\begin{theorem} \label{tm}
Let $(\M,g)$ be a stationary Lorentzian manifold endowed with a
complete timelike Killing vector field $K$ and a complete (smooth,
spacelike) Cauchy hypersurface $S$. Denoting by $\Psi\colon\R\times
\M\to\M$ the flow of $K$, let $P$ and $Q$ be two immersed, disjoint,
connected, closed submanifolds of  $\mathcal M$ which satisfy either
condition $(H_1)$ or $(H_2)$. Then, there exists at least one normal
geodesic joining $P$ to $Q$ in $\mathcal M$.
\end{theorem}

The rest of the paper is organized as follows:
in Section \ref{secfermat} we introduce the Fermat metrics and some basic notions
about semi--Riemannian submersions, while in Section \ref{s2}
we prove Theorem \ref{tm} and discuss some
multiplicity results for geodesics connecting two submanifolds if some suitable assumptions are added to those
in Theorem \ref{tm}.

%-----------------------------------------------------------------

\section{\bf Fermat metrics and submersions }\label{secfermat}
\def\theequation{2.\arabic{equation}}\makeatother
\setcounter{equation}{0}
Before proving our main result, we need some notions  from Finsler geometry  and some basic properties of a
submersion in relation with  stationary spacetimes.
In particular, we recall the Fermat metrics of a standard stationary spacetime,
as introduced in  \cite{cym}.

\begin{definition}{\rm
A {\em Finsler manifold} is a couple $(M,F)$ such that $M$
is a smooth finite dimensional manifold and $F: TM\to[0,+\infty)$ is
a {\em Finsler structure} on $M$, i.e. a function such that
\begin{itemize}
\item[{\sl (i)}] it is continuous on $TM$, $C^{\infty}$ on $TM\setminus 0$ and
it vanishes only on the zero section;
\item[{\sl (ii)}] it is fiberwise positively
homogeneous of degree one, i.e. $F(x,\lambda y)=\lambda F(x,y)$
for all $x\in M$, $y\in T_x M$ and $\lambda>0$;
\item[{\sl (iii)}]
it has fiberwise strictly convex square, i.e.  the matrix
$\left(\frac{1}{2}\frac{\partial^2 (F^2)}{\partial
y^i\partial y^j}(x,y)\right)_{i,j}$ is positive definite for all $(x,y)\in TM\setminus 0$.
\end{itemize}}
\end{definition}

If $(M,F)$ is a Finsler manifold, the {\em length} of a piecewise smooth  curve
$\gamma\colon [a,b]\subset\R \to M$ with respect to the Finsler
structure $F$ is defined by
\[
\ell(\gamma)\ =\ \int_a^b F(\gamma(s),\dot\gamma(s))\ d s.
\]
Hence, the {\em distance} between two
arbitrary points $p, q\in M$ is given by
\[
 {\rm dist}(p,q)= \inf_{\gamma\in {\cal P}(p,q)}\ell (\gamma),
\]
where ${\cal P}(p,q)$ is the set of all piecewise smooth curves $\gamma\colon[a,b]\to M$ with
$\gamma(a)=p$ and $\gamma(b)=q$.

Let us point out that, even if the distance function
with respect to a Finsler structure $F$ is non--negative and satisfies the triangle
inequality, it is not symmetric as, in general, $F$ is non--reversible. Thus,
one has to distinguish between the notions of forward and backward metric balls,
Cauchy sequences and completeness (see \cite[\S 6.2]{bcs} for more details).
Anyway, the topologies generated by the forward and the
backward metric balls coincide with the underlying mani\-fold
topology and a suitable version of the Hopf--Rinow Theorem holds
(see \cite[Theorem 6.6.1]{bcs}).

\begin{theorem}[{\bf Finslerian Hopf--Rinow Theorem}]\label{oph}
Taking a Finsler manifold $(M,F)$, the following statements are equivalent:
\begin{itemize}
\item[{\sl (i)}] the associated Finsler metric is forward (or backward) complete;
\item[{\sl (ii)}] the closed and forward (or backward) bounded subsets of $M$
are compact.
\end{itemize}
Moreover, if (i) or (ii) holds, then any pair of points in $M$ is connected by
a geodesic minimizing the Finslerian distance.
\end{theorem}

A Finsler metric on $M$ is said of {\em Randers type} if
\be\label{finsran}
F(x,y)\ =\ \sqrt{h(x)[y,y]}+\omega(x)[y],
\ee
where $h$ is a Riemannian metric on $M$
and $\omega$ is a one--form such that
\[
\|\omega\|_x\ <\ 1,\quad \hbox{where}\;\; \|\omega\|_x= \sup_{y\in T_x M\setminus \{0\}}\frac{|\omega(x)[y]|}{\sqrt{h(x)[y,y]}}.
\]

Now, let $(\M=\M_0\times\R, g)$ be a standard
stationary Lorentzian manifold, with $g$ as in (\ref{stand}).
If $z(s) = (x(s),t(s))\in{\mathcal M}$, $s\in [0,1]$, is a piecewise smooth future--pointing
or past--pointing lightlike curve, then it satisfies
\be\label{light}
g_0(x)[\dot x,\dot x]+2g_0(x)[\delta(x),\dot x]\dot t-\beta(x)\dot t^2=0.
\ee
Solving equation (\ref{light}) with respect to $\dot t$ and integrating over
 interval $[0,1]$, we get that the difference  $T_{\pm}(x)$  between
the arrival time $t(1)$ and the starting time $t(0)$ of the lightlike curve
$z$ depends only on $x$ and is given by
\be\label{ti}
T_+(x)= \inte \big(\tig(x)[\delta(x),\dot x]+
\sqrt{(\tig(x)[\delta(x),\dot x])^2+\tig(x)[\dot x,\dot x]} \big)\ d s
\ee
if $z$ is future--pointing and by
\be\label{meno}
T_-(x)= \inte \big(\tig(x)[\delta(x),\dot x]-
\sqrt{(\tig(x)[\delta(x),\dot x])^2+\tig(x)[\dot x,\dot x]} \big) \ d s
\ee
if it is past--pointing.  Here, $\tig$ denotes the conformal metric $g_0/\beta$.

\begin{definition}
{\rm The {\em Fermat metrics} associated to $(\M,g)$ are the
Randers metrics $F_+$ and $F_-$ on $\M_0$ respectively given by
\be\label{Fermatmetric}
\begin{split}
F_+(x,y) =& \tig(x)[\delta(x),y]+
\sqrt{(\tig(x)[\delta(x),y])^2+\tig(x)[y,y]}\\
F_-(x,y)=&-\tig(x)[\delta(x),y]+
\sqrt{(\tig(x)[\delta(x),y])^2+\tig(x)[y,y]}
\end{split}
\ee
for every $(x,y)\in T\mathcal M_0$, where the associated Riemannian metric
$h$ in \eqref{finsran} is given by
\[
h(x)[y,y]=(\tig(x)[\delta(x),y])^2+\tig(x)[y,y]
\]
and
\begin{align*}
\omega_+(x)&=\tig(x)[\delta(x),\dot x]&\text{is the one--form related to $F_+$,}\\
\omega_-(x)&=-\tig(x)[\delta(x),\dot x]&\text{is the one--form related to $F_-$.}
\end{align*}
}
\end{definition}
Thus, according to \eqref{ti}-\eqref{Fermatmetric}, if  $z = (x,t)$ is a lightlike curve in $\M$
we have
\[
\Delta_z = t(1) - t(0) = T_{\pm}(x)\ =\ \pm\ \int_0^1 F_\pm(x(s),\dot x(s))\ d s
\]
with $+$ if $z$ is future--pointing, resp. $-$ if $z$ is past--pointing; hence,
$T_\pm(x)$ is $\pm$ the length of the spatial projection $x$ with respect to
the Fermat metric $F_\pm$.

Observe that $F_-$ can be obtained by $F_+$ reversing the sign of $\delta$.
Moreover, if $F_+$ is forward (resp. backward) complete, $F_-$ is backward
(resp. forward) complete and vice versa.

Let us recall the following proposition (cf. \cite[Theorem 4.8]{cym}):

\begin{proposition}\label{burg}
If $(\M = \M_0 \times \R,g)$ is a standard stationary Lorentzian manifold and $\bar{t}\in\R$, then:
\begin{itemize}
\item[$(1)$] if the Fermat  metrics in (\ref{Fermatmetric})
are forward or backward complete on ${\mathcal M}_0$, then
$({\mathcal M},g)$ is globally hyperbolic;
\item[$(2)$]
if $({\mathcal M},g)$ is globally hyperbolic with Cauchy hypersurface
$S={\mathcal M}_0\times \{\bar{t}\}$, then both $F_+$ and $F_-$ in (\ref{Fermatmetric})
are forward and backward complete on ${\mathcal M}_0$.
\end{itemize}
\end{proposition}

Moreover, the following result holds (cf. \cite[Theorem 4.4]{CaJaS09}):

\begin{proposition}
Let $(\M = \M_0 \times \R,g)$ be a standard stationary spacetime and fix $\bar{t}\in\R$.
If $F_+$, or equivalently $F_-$, in (\ref{Fermatmetric}) is forward and backward complete on ${\mathcal M}_0$, then
$S={\mathcal M}_0\times \{\bar{t}\}$ is a Cauchy hypersurface.
\end{proposition}

Hence, the quoted result in \cite{cfs} can be stated as follows: a
standard stationary Lorentzian manifold $(S\times\R,g)$ with
complete Riemannian component $(S,g_0)$ and forward and backward
complete Fermat metric $(S,F_+)$ is geodesically connected.

Furthermore, observe that by \cite[Corollary 3.4]{s} sufficient conditions for the global hyperbolicity of
a standard stationary spacetime $\M = \M_0 \times \R$ are the completeness of
the Riemannian part $({\mathcal M}_0,g_0)$ and some
growth assumptions on the coefficients $\delta, \beta$ of the metric (\ref{stand}).
As shown in \cite{bcf}, such conditions are optimal to get geodesic connectedness on standard stationary spacetimes. Indeed, in \cite{bcf}
it is furnished an example where variational techniques, commonly employed to prove geodesics connectedness of standard stationary spacetimes, fail
and the presumable lack of connectedness by geodesics
can be explained by a geometric viewpoint by the fact that the associated Fermat metrics are
not forward and backward complete.

We conclude this section introducing the basic notions on
semi--Riemannian submersions needed in the proof of Theorem~\ref{tm} in the hypothesis $(H_2)$ (e.g., cf. \cite{o}).

\begin{definition}
{\rm Let $(M,g_1)$ and $(B,h_1)$ be two semi--Riemannian manifolds.
A {\sl semi--Rie\-mann\-ian submersion} between $M$ and $B$ is a smooth map $\pi\colon M\to B$ such that for any $p\in M$:
\begin{itemize}
\item[{\sl (i)}]
 its differential $\mathrm{d} \pi(p)\colon T_pM\to T_{\pi(p)}B$ is surjective;
\item[{\sl (ii)}]
the fiber $\pi^{-1}(\pi(p))$ is a non-degenerate submanifold  of $M$;
\item[{\sl (iii)}]
$\mathrm d\pi(p)\colon \mathcal H T_p M\to T_{\pi(p)} B$ is an isometry.
\end{itemize}
Here, $\mathcal H T_p M$ is the {\em horizontal} subspace of $T_pM$,
i.e. the orthogonal subspace  to the {\em vertical} one $\mathcal V T_pM=\mathrm{Ker}(\mathrm d\pi(p))$.}
\end{definition}
Giving a $C^1$ curve $\gamma:[a,b]\to B$, a \emph{horizontal lift}
of $\gamma$ is a curve
$\alpha:[a,c]\subset [a,b]\to M$
such that $\pi\circ\alpha=\gamma$
and $\dot\alpha(s)$ is horizontal, i.e.
$\dot\alpha(s)\in \mathcal H T_{\alpha(s)} M$ for all $s\in [a,c]$.

If  $(\mathcal M=\M_0\times \R, g)$ is a standard stationary spacetime, then the canonical
 projection $\pi_{\M_0}$ on $\M_0$ is a Lorentzian submersion between $(\mathcal M, g)$ and $(\M_0,h_1)$, where $h_1$
is the Riemannian metric defined as
\begin{equation}\label{h1}
h_1(x)[v,v]=g_0(x)[v,v]+\tfrac{1}{\beta(x)}(g_0(x)[\delta(x),v])^2.
\end{equation}
In fact, writing the metric $g$ as
\[\begin{split}
g(x,t)[(v,\tau),(v,\tau)]\ =\ &g_0(x)[v,v]+\tfrac{1}{\beta(x)}(g_0(x)[\delta(x),v])^2\\
&-\left (\tfrac{1}{\sqrt{\beta(x)}}g_0(x)[\delta(x),v]-\sqrt{\beta(x)}\tau\right )^2
\end{split}
\]
and considering $\mathcal H T_{(x,t)}\mathcal M$, the orthogonal subspace to the one--dimensional subspace
$[\partial_t|_{(x,t)}]$
for all $(x,t)\in\M_0\times\R$, the map $\mathrm{d} \pi_{\M_0}(x,t): \mathcal H T_{(x,t)}\mathcal M\to T_x\mathcal M_0$
is an isometry with respect to the restriction of $g(x,t)$ to $ \mathcal H T_{(x,t)} \mathcal M$ and $h_1(x)$.

%-----------------------------------------------------------
\section{\bf Proof of Theorem~\ref{tm}}\label{s2}
Throughout this section,
$(\M,g)$ is a stationary spacetime endowed with a
complete timelike Killing vector field $K$
and $P$ and $Q$ are two disjoint connected closed immersed submanifolds of $\mathcal M$.

Firstly, let us point out that $\M$ can be equipped with a Riemannian metric defined as follows:
\[
g_R(p)[v_1, v_2]\ =\ g(p)[v_1, v_2]\ -\ 2\ \frac{g(p)[v_1, K(p)]\ g(p)[ v_2, K(p)]}{ g(p)[K(p) , K(p) ]}
\]
for all $p \in \mathcal M$, $v_1 , v_2 \in T_p \mathcal M$.

Furthermore,  by standard arguments
it can be proved that normal geodesics joining $P$ to $Q$ are the critical points of the functional
\be\label{fu}
f(z)\ =\ \frac 12\int_0^1 g(z)[\dot z,\dot z]\ d s
\ee
defined on the Hilbert manifold $\Omega(P,Q)$ of the $H^1-$curves connecting $P$ to $Q$, that is
\[
\begin{split}
\Omega(P,Q)\ =\ \big\{z:[0,1]\rightarrow \mathcal M :\, &z\ \text{ is absolutely continuous,}\\
& z(0) \in P,\, z(1) \in Q,\, \int_0^1 g_R(z)[\dot z,\dot z] \ d s < +\infty\big\},
\end{split}
\]
so for each $z\in\Omega(P,Q)$ the tangent space $T_z\Omega(P,Q)$ is given by the $H^1$--vector fields
$\zeta\colon[0,1]\to T\mathcal M$ along $z$ such that $\zeta(0)\in T_{z(0)}P$ and $\zeta(1)\in T_{z(1)}Q$.

In our setting, geodesics satisfy the conservation law \eqref{vinc},
hence the critical points of the functional $f$ belong to the subset  
\[
\Omega_K(P,Q) = \left\{ z \in \Omega (P,Q) :\ \exists\, C_z \in \R \ \hbox{s.t.}\
g(z)[\dot z, K(z)]= C_z\ \text{a.e. on $[0,1]$}\right\}.
\]
We point out that $\Omega_K(P,Q)$ may be empty. In order to guarantee $\Omega_K(P,Q) \ne \emptyset$,
it is enough to assume that $K$ is complete (see \cite[Lemma 5.7]{gp} and \cite[Proposition 3.6]{cfs})
or to consider two submanifolds $P$ and $Q$ causally related (see \cite[Appendix A]{bgs}).

Furthermore, it can be proved not only that $\Omega_K(P,Q)$ is a smooth submanifold of $\Omega(P,Q)$
(see \cite[Proposition 3.1]{bgs}) but also that a good variational principle holds on it:

\begin{theorem} \label{t1}
Let $f$ be as in (\ref{fu}).
A curve $z\in \Omega(P,Q)$ is a critical point of $f$ on $\Omega(P,Q)$ (hence, a normal geodesic connecting $P$ to $Q$)
if and only if it is a critical point of $f$ on $\Omega_K(P,Q)$.
\end{theorem}
\begin{proof}
Obviously, from \eqref{vinc} if $z\in
\Omega(P,Q)$ is a critical point of $f$ on $\Omega(P,Q)$, then it is
a critical point of $f$ on $\Omega_K(P,Q)$.
Now, assume that $z\in
\Omega_K(P,Q)$ is a critical point of $f$ on $\Omega_K(P,Q)$, i.e.
$\de f(z)[\zeta] = 0$ for all $\zeta \in T_z\Omega_K(P,Q)$.
The proof in the particular case $P=\{p\}$ and $Q=\{q\}$, with $p$, $q
\in \M$, is contained in \cite{gp} and is based on the fact that,
for any $z\in \Omega_K(p,q)$, $T_z\Omega(p,q)$ splits into the direct
sum of $T_z\Omega_K(p,q)$ and the space $\mathcal W_z$ of the
vector fields in $T_z\Omega(p,q)$ which are pointwise collinear to
$K(z)$. A vector field $\zeta\in T_z\Omega(p,q)$ belongs to
$\mathcal W_z$ if and only if there exists a function $\mu\in
H^1_0([0,1],\R)$ such that $\zeta=\mu K(z)$. Thus, by
straightforward computations, we get that $\mathrm d f(z)[\zeta]=0$
for all $\zeta\in\mathcal W_z$.
In \cite{bgs} this result is
extended to the more general case of two spacelike submanifolds $P$
and $Q$. Also in this setting, it is \be\label{split1}
T_z\Omega(P,Q)\ =\ T_z\Omega_K(P,Q) \oplus \mathcal W_z, \ee where
\begin{multline*}
T_z \Omega_K(P,Q)=\big\{\zeta\in T_z\Omega(P,Q):\ \exists\, C_{\zeta}\in\R \ \hbox{s.t.}\\
g(z)[\nabla_s\zeta,K(z)]-g(z)[\zeta,\nabla_sK(z)]=C_{\zeta}\ \text{ a.e. on $[0,1]$}\big\}
\end{multline*}
and again
\[
\mathcal W_z = \big\{\zeta \in T_z \Omega(P,Q):\ \zeta= \mu K(z)\ \hbox{with}\ \mu\in  H^1_0([0,1],\R)\big\}.
\]
In general, if $P$ and $Q$ are not spacelike submanifolds, the direct sum \eqref{split1} does not hold.
In fact, if $K(z(0)) \in T_{z(0)}P$ and $K(z(1)) \in T_{z(1)}Q$, then any vector field $\zeta = \mu K(z)$,
with
\[
\mu(s)\ =\ \mu_0+\int_0^s\frac{C}{g(z)[K(z),K(z)]}\ d \tau,\qquad \mu_0, C\in \R,
\]
belongs to $T_z\Omega_K(P,Q)\cap\mathcal W_z$.
Anyway, since the proof of
such a decomposition does not rely on the boundary conditions that
$\zeta$ and $\tilde \zeta$  have to satisfy (see the proof of
\cite[Proposition 3.3]{bgs}), taking any
$\zeta\in T_z\Omega(P,Q)$, a vector field $\tilde\zeta\in
T_z\Omega_K(P,Q)$ and an $H^1_0$--function $\mu$ exist such that
\be\label{dec} \zeta=\tilde\zeta + \mu K(z). \ee 
On the other hand, reasoning as in
\cite[Proposition 2.2]{bgs}, it can be proved that $\de f(z)[\zeta]
= 0$ for all $\zeta \in {\cal W}_z$; hence, \eqref{dec} implies $\de
f(z)[\zeta] = 0$ for all $\zeta \in T_z\Omega(P,Q)$.
\end{proof}

From now on, let us assume that the hypotheses of Theorem~\ref{tm} hold.
Hence, as already remarked in Section \ref{s1}, as $S$ is
a smooth complete Cauchy hypersurface, $\M$ is a standard stationary spacetime which globally splits as
$\mathcal M=S\times\R$ and the metric $g$ is as in (\ref{stand})
with the identifications in (\ref{iden}).
Hence, we have
\[
g(z)[\dot z,K(z)] =
g_0(x)[\delta(x),\dot x] - \beta(x)\dot t
\]
for any absolutely continuous curve $z=(x,t):[0,1]\rightarrow
\mathcal M$. Thus, if $z\in\Omega_K(P,Q)$, a constant $C_z$ exists
such that \be\label{dott}
 \dot t=\tilde g_0(x)[\delta(x),\dot x]-\frac{C_z}{\beta(x)}
\ee
(as in Section \ref{secfermat}, we set $\tig=g_0/\beta$).
Integrating both hand sides of (\ref{dott}) in $[0,1]$, we get
\be \label{cz}
C_z= \left(\int_0^1\tilde g_0(x)[\delta(x),\dot x] d s-\Delta_z\right)
\left(\int_0^1\frac{1}{\beta(x)} d s\right)^{-1},
\ee
with $\Delta_z=t(1)-t(0)$.  Now, replacing (\ref{stand})
in (\ref{fu}) and substituting (\ref{cz}) in (\ref{dott}),  we can express
the restriction of $f$ to $\Omega_K(P,Q)$, denoted by $\J$, as a functional
depending only on the $x$ component of the curve $z\in\Omega_K(P,Q)$
and on $\Delta_z$ (for the first claim of this variational principle, see \cite{gm}):
\be\label{J}
\begin{split}
\J(z)= &\frac 12\int_0^1g_0(x)[\dot x,\dot x] \ d s
+\frac 12\int_0^1\tilde g_0(x)[\delta(x),\dot x]\ g_0(x)[\delta(x),\dot x] \ d s\\
&-\frac 12\left(\int_0^1 \tilde g_0(x)[\delta(x),\dot x] \ d s -
\Delta_z\right)^2 \left(\int_0^1 \frac1{\beta(x)}\ d s\right)^{-1}.
\end{split}
\ee

We recall  that a $C^1$ functional $J\colon \Omega\to\R$, defined on
a Hilbert manifold $\Omega$, satisfies the {\em Palais--Smale condition} if
each sequence $(z_n)_n\subset \Omega$, such that $(J(z_n))_n$ is bounded
and $\de J(z_n)\to 0$ admits a converging subsequence.

A classical existence theorem for critical points of a functional defined on a Hilbert manifold is the following.

\begin{theorem}\label{abstract}
If $\Omega$ is a Hilbert manifold and $J: \Omega \to \R$ is a $C^1$ functional which  satisfies the Palais--Smale condition,
is bounded from below and has a complete non--empty sublevel, then it attains its infimum.
\end{theorem}

As our aim is applying the previous abstract theorem to $\J$ in $\Omega_K(P,Q)$, firstly we state a pair of remarks useful in the proof
of the boundedness from below and of the Palais--Smale condition for such a functional.

\begin{remark}\label{due}
{\rm Note that under assumption $(H_1)$ of Theorem \ref{tm}, in the
splitting $\M = S \times \R$ (recall \eqref{iden}),  from \eqref{sQ}
and the compactness of $P$ we have:
\[
|\Delta_z| =|t(1)-t(0)|\leq D_Q+D_P,
\]
(here, $\sup|s_P(P)|=D_P<+\infty$); hence $\Delta_z$ is bounded on
$\Omega_K(P,Q)$.}
\end{remark}
\begin{remark}\label{waw}
{\rm In Section \ref{secfermat}, associated to any piecewise smooth
lightlike curve $s\mapsto z(s)=(x(s),t(s))\in S\times\R$ we  have
introduced the quantities $T_\pm(x)$, depending only on $x$ (see
\eqref{ti} and \eqref{meno}), so that $\Delta_z = T_+(x)$ if $z$ is
future--pointing, $\Delta_z = T_-(x)$ if $z$ is past--pointing.

 On the other hand, as in \cite{fgm}, we can consider \be\label{pp} z =
(x,t) \in \Omega_K(P,Q) \quad \hbox{such that $\J(z) = 0$,} \ee then
it has to be
\[
\begin{split}
\Delta_z\ = &\int_0^1 \tilde g_0(x)[\delta(x),\dot x]\ d s\\
&+
\sqrt{\left(\|\dot x\|^2 +
\int_0^1\tilde g_0(x)[\delta(x),\dot x]\ g_0(x)[\delta(x),\dot x]\ d s\right)
 \int_0^1 \frac1{\beta(x)}\ d s}
\end{split}
\]
if $\Delta_z > 0$, while
\[
\begin{split}
\Delta_z\ = &\int_0^1 \tilde g_0(x)[\delta(x),\dot x]\ d s\\
&-
\sqrt{\left(\|\dot x\|^2 +
\int_0^1\tilde g_0(x)[\delta(x),\dot x]\ g_0(x)[\delta(x),\dot x] \ d s\right)
 \int_0^1 \frac1{\beta(x)} \ d s}
\end{split}
\]
if $\Delta_z < 0$, where we have set
\[\|\dot x\|^ 2= \int_0^1g_0(x)[\dot x,\dot x] \ d s.\]
Thus, depending only on $x$, we can define
\be\label{titilde}
\begin{split}
&\widetilde T_+(x)\ =\ \int_0^1 \tilde g_0(x)[\delta(x),\dot x] \ d s\\
&\qquad +
\sqrt{\left(\|\dot x\|^2 +
\int_0^1\tilde g_0(x)[\delta(x),\dot x]\ g_0(x)[\delta(x),\dot x] \ d s\right)
 \int_0^1 \frac1{\beta(x)} \ d s},
\end{split}
\ee
\be\label{menotilde}
\begin{split}
&\widetilde T_-(x)\ =\ \int_0^1 \tilde g_0(x)[\delta(x),\dot x] \ d s\\
&\qquad -
\sqrt{\left(\|\dot x\|^2 +
\int_0^1\tilde g_0(x)[\delta(x),\dot x]\ g_0(x)[\delta(x),\dot x] \ d s\right)
 \int_0^1 \frac1{\beta(x)} \ d s},
\end{split}
\ee so that, if \eqref{pp} holds, we have $\Delta_z = \widetilde
T_\pm(x)$ according to the sign of $\Delta_z$.

 Observe that, for any
$z = (x,t) \in \Omega_K(P,Q)$ (non-necessarily lightlike), by
comparing \eqref{ti} with \eqref{titilde} and \eqref{meno} with
\eqref{menotilde}, the definition of $\tilde g_0$ and the
Cauchy--Schwarz inequality imply \be\label{stimo} T_+(x) \le
\widetilde T_+(x),\quad \widetilde T_-(x) \le T_-(x). \ee}
\end{remark}

\begin{theorem}\label{bbps}
Under the hypothesis $(H_1)$ of Theorem~\ref{tm}, the functional $\J$ is bounded from below, has complete sublevels and
satisfies the  Palais--Smale condition in $\Omega_K(P,Q)$.
\end{theorem}
\begin{proof}
Let us divide the proof in three steps.
Firstly, we claim that, taking a sequence $(z_n)_n\subset
\Omega_K(P,Q)$ such that \be\label{star} (\J(z_n))_n \quad \hbox{is
bounded from above}, \ee and considering the splitting $S\times \R$,
so that $z_n=(x_n,t_n)$, we have
\be\label{cl} (\|\dot x_n\|)_n \;
\hbox{is bounded.} \ee In fact, arguing by contradiction, let us
assume that, up to subsequences, \be\label{div} \|\dot x_n\| \
\xrightarrow{n}\ +\infty. \ee From (\ref{J}) and the Cauchy--Schwarz
inequality it follows \be\label{2J}
\begin{split}
&2 \J(z_n)\ \ge\ \|\dot x_n\|^2\\
&\quad - \Delta_{z_n} \left(\Delta_{z_n} - 2 \int_0^1\tilde g_0(x_n)[\delta(x_n),\dot x_n]
\ d s\right) \left( \int_0^1\frac{1}{\beta(x_n)} \ d s\right)^{-1},
\end{split}
\ee where, by Remark \ref{due}, we have \be\label{del}
(|\Delta_{z_n}|)_n\quad \hbox{is bounded.} \ee Hence, we can rule
out both the following possibilities: the fact that $(\Delta_{z_n})_n$ is definitively equal to
$0$ and the existence of a compact subset of $S$ containing all the
supports of the curves $x_n$, $n\in \N$;  otherwise, from
\eqref{div} and \eqref{2J}, it would follow \be\label{div2} \J(x_n)\
\xrightarrow{n}\ +\infty, \ee in contradiction with \eqref{star}.
Thus, let us assume that no compact subset of $S$ contains the images of
all the curves $x_n$
and, up to subsequences, $\Delta_{z_n}>0$, resp. $\Delta_{z_n}<0$, for all $n\in\mathbb N$.
Then, a subsequence exists such that \be\label{timeno}
T_+(x_n)\xrightarrow{n} +\infty \hbox{ if } \Delta_{z_n}>0,\qquad
\hbox{resp. }T_-(x_n)\xrightarrow{n} -\infty \hbox{ if }
\Delta_{z_n}<0 \ee (recall (\ref{ti}) and (\ref{meno})). In fact,
since $S$ is a Cauchy hypersurface, from Proposition~\ref{burg} (2),
the Fermat metrics defined in (\ref{Fermatmetric}) on $S$ are
forward and backward complete, then by the already recalled
Finslerian Hopf--Rinow Theorem (see Theorem \ref{oph}), if
$(T_+(x_n))_n$ is bounded from above (resp. $(T_-(x_n))_n$ is
bounded from below), as the sequence $(x_n(0))_n$ is contained in
the compact subset $\Psi(G_P)$ (see \eqref{proiezionecompatta}), a
compact subset of $S$ must contain all the images of the curves
$x_n$, which is a contradiction.
Moreover, if $\Delta_{z_n}>0$,
resp. $\Delta_{z_n}<0$, according to Remark \ref{waw}, related to
each $z_n$ it can be considered $\widetilde T_+(x_n)$, resp.
$\widetilde T_-(x_n)$, and from \eqref{stimo}, \eqref{timeno} it
follows \be\label{timeno2} \widetilde T_+(x_n)\xrightarrow{n}
+\infty \hbox{ if } \Delta_{z_n}>0,\qquad \hbox{resp. }\widetilde
T_-(x_n)\xrightarrow{n} -\infty \hbox{ if } \Delta_{z_n}<0. \ee But
from \eqref{J} and \eqref{titilde} it follows
\[\begin{split}
2\J(z_n) = &(\widetilde T_\pm(x_n) - \Delta_{z_n})\
\big(\widetilde T_\pm(x_n) + \Delta_{z_n}\\
& - 2 \int_0^1\tilde g_0(x)[\delta(x),\dot x] \ d s\big)\
 \big(\int_0^1 \frac1{\beta(x)} \ d s\big)^{-1},
 \end{split}
\]
with the $+$ sign in $\widetilde T_\pm$ if $\Delta_{z_n}>0$, resp.
the $-$ sign if $\Delta_{z_n}<0$; so, following \cite[Lemma
5.6]{cfs} by (\ref{div}), \eqref{del} and (\ref{timeno2}) we get
\eqref{div2} in contradiction with \eqref{star}. Therefore claim
(\ref{cl}) is proved and, as $S$ is complete with respect to the
metric $g_0$ and the sequence $(x_n(0))_n$ is contained in the
compact set in
\eqref{proiezionecompatta}, all the supports of the curves $x_n$ lie in a compact subset of $S$.
Now, we can prove that $\J$ is bounded from below in $\Omega_K(P,Q)$. In fact,
taking a minimizing sequence $(z_n)_n\subset \Omega_K(P,Q)$ for $\J$, namely
\[
 \lim_{n\rightarrow +\infty}\J(z_n)=\inf_{z\in \Omega_K(P,Q)} \J(z),
\]
we have that \eqref{star}, hence \eqref{cl}, holds. Thus, from \eqref{2J} and
\eqref{del} we have that the sequence $(\J(z_n))_n$ is bounded from below, too, whence
\[
\inf_{z\in \Omega_K(P,Q)} \J(z) > -\infty.
\]
At last, we have to prove that $\J$ satisfies the Palais--Smale
condition. To this aim, let $(z_n)_n\subset \Omega_K(P,Q)$,
$z_n=(x_n,t_n)$ according to the splitting $S\times \R$, be such
that $(\J(z_n))_n$ is bounded and $\de \J(z_n)\to 0$.
Obviously \eqref{star} holds; so as above, the components $x_n$ satisfy
\eqref{cl} and have supports contained in a compact subset of $S$.
Hence, by the Ascoli--Arzel\`a Theorem a uniformly convergent
subsequence of $(x_n)_n$ exists. Furthermore, by (\ref{cz}) also the
sequence $(C_{z_n})_n$ is bounded and, by (\ref{dott}), so it is for $(\|\dot
t_n\|)_n$. As $(t_n(0))_n$ is contained in a compact subset of $\R$,
again by the Ascoli--Arzel\`a Theorem there exists also a
subsequence of $(t_n)_n$ which uniformly converges. Then the
existence of a subsequence converging in $\Omega_K(P,Q)$ and the
completeness of the sublevels of $\J$ can be obtained respectively
as in \cite[Theorem 5.1]{bgs} and in \cite[Proposition 5.2]{bgs}.
\end{proof}
\begin{proof}[Proof of Theorem~\ref{tm}]
Under assumption
$(H_1)$, the existence of a minimum of $\J$ in $\Omega_K(P,Q)$
follows from Theorems \ref{abstract} and \ref{bbps}. Hence, by
Theorem~\ref{t1} such a minimum is a normal geodesic connecting $P$
and $Q$.
In the case $(H_2)$, recalling that the canonical projection
$\pi_S\colon (S\times\R, g)\to (S, h_1)$, where $h_1$ is the metric
defined in \eqref{h1}, is a Lorentzian submersion, we can use the
fact that the horizontal lift of any geodesic in the base of a
semi--Riemannian submersion is a geodesic of the total space (see
\cite[Corollary 7.46]{o}). Since $g_0$ is complete, also $h_1$ is
complete. From a theorem of K. Grove \cite[Theorem 2.6]{g}, at least
one normal geodesic $x:[0,1]\to S$ in $(S,h_1)$ connecting $P_S$ and
$Q_S$ exists. Hence, a horizontal lift of such a geodesic provides a
normal geodesic of $(\M, g)$ connecting $P$ to $Q$ (observe that the
$t$ component of its horizontal lift is given by $t(s)= t_0 +
\int_0^s\frac 1{\beta(x)}g_0(x)[\delta(x),\dot x]\ d \tau$).
\end{proof}

Let us point out that under assumption $(H_2)$ of Theorem~\ref{tm},
the normal geodesic connecting $P$ and $Q$, being horizontal, is spacelike. Moreover,
changing the initial point of the geodesic on the fiber we obtain infinitely many
spacelike normal geodesics connecting $P$ to $Q$ which all project on the same
geodesic on $(S,h_1)$.

A more interesting multiplicity result can be obtained minimizing
the energy functional of the Riemannian manifold $(S,h_1)$ on
homotopy classes of curves from $P_S$ to $Q_S$ or assuming
that $S$ is not contractible and $P_S,Q_S$ are contractible in $S$.
In this last case the Ljusternik--Schnirelmann category of $\Omega(P_S,Q_S)$ is
infinite (cf. \cite{cc} and \cite{fh}) and then infinitely many
normal geodesics in $(S, h_1)$ connecting $P_S$ and $Q_S$ exist (see \cite[Theorem 2.6]{g});
therefore, there exist infinitely many spacelike normal geodesics
in $(\M,g)$ from $P$ to $Q$, having different projections on $S$ (up to be
the iterates of a closed prime geodesic of $(S, h_1)$ crossing orthogonally $P_S$ and $Q_S$).

An analogous multiplicity result can be also obtained under the assumption $(H_1)$. Indeed, if
the Killing vector field $K$ is complete the manifold $\Omega_K(P,Q)$ is homotopically
equivalent to $\Omega(P,Q)$ (see \cite[Proposition 5.5]{bgs})  and we can apply again
Ljusternik--Schnirelmann Theory if suitable hypotheses on $\M$,
$P$ and $Q$ imply that the category of $\Omega(P,Q)$ is non--trivial.

%----------------------------------------------------------


\small

\begin{thebibliography}{777}

\bibitem{bcf} {R. Bartolo, A.M. Candela and J.L. Flores,}
{\em Geodesic connectedness of stationary spacetimes
with optimal growth}, J. Geom. Phys. {\bf 56} (2006), 2025-2038.

\bibitem{bgs} {  R. Bartolo, A. Germinario and M. S\'anchez,}
{\em Orthogonal trajectories on stationary spacetimes under intrinsic assumptions},
Topol. Methods Nonlinear Anal. {\bf 24} (2004), 239-268.

\bibitem{bgs1}
{R. Bartolo, A. Germinario and M. S\'anchez,}
{\em Trajectories connecting two submanifolds on non--complete Lorentzian manifolds},
Electron. J. Differential Equations {\bf 10} (2004), 20 pp.

\bibitem{bcs}
D. Bao, S.S. Chern and Z. Shen,  An Introduction to {R}iemann-{F}insler Geometry.
Graduate Texts in Mathematics, Springer--Verlag, New York, 2000.

\bibitem{bee}
 J.K. Beem, P.E. Ehrlich and K.L. Easley,
Global Lorentzian Geometry, 2nd Ed.,
Pure App. Math. \textbf{202}, Marcel Dekker, New York, 1996.

\bibitem{c} {A.M. Candela,}
{\em Normal geodesics in static spacetimes with critical asymptotic
behavior}, Nonlinear Anal. TMA {\bf 63} (2005), 357-367.

\bibitem{cfs} {A.M. Candela, J.L. Flores and M. S\'anchez,}
{\em Global hyperboliticity and Palais--Smale condition
for action functionals in stationary spacetimes}, Adv. Math. {\bf 218} (2008), 515-536.

\bibitem{cms}
A.M. Candela, A. Masiello and A. Salvatore,
{\em Existence and multiplicity of normal geodesics in {L}orentzian manifolds},
J. Geom. Anal. {\bf 10} (2000), 623-651.

\bibitem{cs} A.M. Candela and A. Salvatore, {\em Light rays joining two submanifolds in space-times},
J. Geom. Phys. {\bf 22} (1997), 281-297.

\bibitem{cs1} A.M. Candela and A. Salvatore,
{\em Normal geodesics in stationary Lorentzian manifolds with unbounded
coefficients}, {J. Geom. Phys.} {\bf 44} (2002), 171-195.

\bibitem{cc}
A. Canino,
{\em On $p$-convex sets and geodesics},
J. Differential Equations {\bf 75} (1988), 118-157.

\bibitem{cym} E. Caponio, M.A. Javaloyes and A. Masiello, {\em On the energy functional 
on Finsler manifolds and applications to stationary spacetimes}, arXiv: math/0702323v3 [math.DG] (2008).

\bibitem{CaJaPi09} {E. Caponio, M.A. Javaloyes and P. Piccione,} {\em Maslov index in semi-Riemannian  
submersions},  
Ann. Global Anal. Geom. {\bf 38} (2010), 57-75.

\bibitem{CaJaS09} {E. Caponio, M.A. Javaloyes and M. S{\'a}nchez,} {\em On the interplay between Lorentzian causality and Finsler metrics of Randers type},
arXiv:0903.3501v1 [math.DG] (2009).

\bibitem{fh}
  E. Fadell and S. Husseini,
  {\em Category of loop spaces of open subsets in {E}uclidean space},
  Nonlinear Anal. {\bf 17} (1991), 1153-1161.

\bibitem{fgm}
D. Fortunato, F. Giannoni and A. Masiello,
{\em A {F}ermat principle for stationary space-times and applications
to light rays}, J. Geom. Phys. {\bf 15} (1995), 159-188.

\bibitem{gm}
F. Giannoni and A. Masiello,
{\em On the existence of geodesics on stationary {L}orentz manifolds with
convex boundary},
J. Funct. Anal. {\bf 101} (1991), 340-369.

\bibitem{gp}
F. Giannoni and P. Piccione, {\em An intrinsic approach to the geodesical
connectedness of stationary Lorentzian manifolds}, Comm. Anal. Geom. {\bf 1} (1999), 157-197.

\bibitem{g}
K. Grove, {\em Condition {$(C)$} for the energy integral on certain path spaces and applications to the theory of geodesics},
J. Differential Geometry \textbf{8} (1973), 207-223.

\bibitem{JavS08}
{M.A. Javaloyes and S. S{\'a}nchez},
{\em A note on the existence of standard splittings for conformally stationary
spacetimes}, Classical Quantum Gravity \textbf{25} (2008), 168001  (7 pp).

\bibitem{o}
{B. O'Neill}, Semi--Riemannian Geometry with Applications to Relativity, Academic Press, New York--London, 1983.

\bibitem{s} M. S\'anchez,
{\em Some remarks on causality theory and variational methods in
Lorentzian manifolds}, Conf. Semin. Mat. Univ. Bari {\bf 265} (1997).

\bibitem{s1}
{M. S\'anchez}, {\em Lorentzian manifolds admitting a Killing vector field}, Nonlinear Anal. TMA {\bf 1} (1997), 643-654.

\end{thebibliography}
\end{document}